\documentclass[a4paper,10pt]{amsart}
\usepackage{geometry,float}
\usepackage{bbm,amsmath,amssymb,amsfonts,amsthm}
\usepackage{enumerate}
\usepackage{a4}
\usepackage{url}

\newcommand{\NN}{\mathbbm{N}}    %

\newcommand{\ZZ}{\mathbbm{Z}}    %
\newcommand{\QQ}{\mathbbm{Q}}    %
\newcommand{\RR}{\mathbbm{R}}    %

\newcommand{\DD}{\mathcal{D}}
\newcommand{\GG}{\mathcal{G}}
\newcommand{\lieg}{\mathfrak{g}}

\newcommand{\CC}{\mathbbm{C}}   %

\renewcommand{\phi}{\varphi}
\renewcommand{\epsilon}{\varepsilon}

\DeclareMathOperator{\diam}{diam} 
 
\DeclareMathOperator{\Fix}{Fix}

\DeclareMathOperator{\Addim}{Addim}
\DeclareMathOperator{\Stab}{Stab}

\newcommand{\OO}{\mathcal O}
\newcommand{\UU}{\mathcal U}
\newcommand{\mA}{\mathcal A}
\DeclareMathOperator{\vol}{vol}
\DeclareMathOperator{\QH}{QH}
\DeclareMathOperator{\Lie}{Lie}

\DeclareMathOperator{\Hom}{Hom}

\DeclareMathOperator{\Isom}{Isom}

\newcommand{\mcL}{{\mathcal L}}

\DeclareMathOperator{\im}{im}

\DeclareMathOperator{\diag}{diag}

\DeclareMathOperator{\Aut}{Aut}

\DeclareMathOperator{\GL}{GL}
\DeclareMathOperator{\Sym}{Sym}

\DeclareMathOperator{\SL}{SL}

\DeclareMathOperator{\Ad}{Ad}

\newtheorem{proposition}{Proposition}[section]
\newtheorem{theorem}[proposition]{Theorem}
\newtheorem{lemma}[proposition]{Lemma}
\newtheorem{corollary}[proposition]{Corollary}

\theoremstyle{remark}
\newtheorem{example}[proposition]{Example}

\newtheorem{remark}[proposition]{Remark}

\theoremstyle{definition}
\newtheorem{definition}[proposition]{Definition}

\begin{document}

\title{Embeddings of algebraic groups in Kac--Moody groups}
\author{Guntram Hainke}
\email{guntram.hainke@math.uni-giessen.de}

\maketitle
\begin{abstract} Let $k_1,k_2$ be two fields of characteristic 0.
Let $G_1$ be a split semisimple algebraic group over $k_1$, $G_2$ a split Kac--Moody group over $k_2$ and $\varphi\colon G_1(k_1)\to G_2(k_2)$ an abstract embedding. We show that $\im \varphi$ is a bounded subgroup whenever $k_1$ is an algebraic extension of $\QQ$, while there are embeddings with unbounded image if $k_1$ has infinite transcendence degree over $\QQ$. 
\end{abstract}

\section{Introduction}
Let $k$ be a field, $\mathcal D$ a Kac--Moody root datum and $G:=\mathcal G_\mathcal{D}(k)$ a Kac--Moody group as constructed by Tits \cite{Tits87}. The group $G$ is given by a presentation which generalizes Steinberg's presentation of the $k$-rational points of a split semisimple algebraic $k$-group.\\
The isomorphism problem for Kac--Moody groups asks to describe all possible isomorphisms $\varphi\colon G \to G'$, where $G':=\mathcal G_{\mathcal D'}(k')$ is another Kac--Moody group over a possibly different field $k'$. 
It was investigated in \cite{MR1155464}, \cite{PECAbstract}, \cite{MR2261755} and \cite{Hainke} under increasingly weaker hypotheses. \\
For two parabolic subgroups $P_1,P_2 \leq G$ of opposite sign, the intersection $P_1\cap P_2$ is isomorphic to the $k$-rational points of an algebraic $k$-group, and any subgroup contained in such an intersection is called a bounded subgroup. Now an isomorphism $\varphi\colon G \to G'$ induces embeddings of these algebraic groups. Under the condition that  $\varphi$ maps bounded subgroups of $G$ to bounded subgroups of $G'$ (which is automatic over finite fields), Caprace--M\"uhlherr \cite{MR2261755} could show that such an isomorphism is standard, i.e. the product of a field automorphism, a graph automorphism, a diagonal automorphism and a sign automorphism. \\
In this paper, which is based on a part of the author's Ph.D. thesis \cite{Hainke}, we address the related problem of determining when the image of an embedding $\varphi\colon G_1(k_1) \to G_2(k_2)$ of a semisimple algebraic group $G_1(k_1)$ into a Kac--Moody $G_2(k_2)$ has bounded image. We show that if $k$ has infinite transcendence degree over its prime field and $G_1$ is isotropic over $k$, there exist embeddings with unbounded image, whereas the image is necessarily bounded if $G_1$ is split over $k_1$ and $k_1$ is an algebraic extension of the rational numbers. \\

{\bf Acknowledgement.} I am grateful to Pierre-Emmanuel Caprace for several useful discussions on the subject of this article.\\

\section{Preliminaries} 
\subsection{Generalities about bounded actions}
\index{bounded!generation}
Let $G$ be a group and $(U_i)_{i \in I}$ a family of subgroups of $G$. Then $G$ is said to be {\bf boundedly generated} by $(U_i)_{i \in I}$ if there exists $n \in \NN$ such that each $g \in G$ is a product $g=g_{i_1} \ldots g_{i_n}$ of at most $n$ elements $g_{i_j} \in \cup_{i \in I} U_i$.\\
Let $X$ be a complete CAT(0)-space (cf. \cite{BridsonHaefliger}) and let $G$ be a group which acts by isometries on $X$. Then the $G$-action is called {\bf bounded} if there is a bounded
$G$-orbit in $X$, i.e. there is some $x \in X$ such that $\diam G \cdot x <\infty.$ By the Bruhat--Tits fixed point theorem (\cite[Corollary II.2.8]{BridsonHaefliger}), this is equivalent to the fact that $G$ fixes a point of $X$. \\

\begin{lemma}\label{boundedgenerationlemma}\label{NormalizesCommonFixedPoint}
Let $X$ be a complete CAT(0)-space.
\begin{enumerate}
\item Let $G\leq \Isom(X)$ be a group which is boundedly generated by a finite number of subgroups $U_1, \ldots, U_n$. Then the $G$-action is bounded if and only if every $U_i$-action is bounded. 
\item Let $N,U \leq \Isom(X)$ be such that the $N$- and the $U$-action are bounded. If $N$ normalizes $U$, then $N$ and $U$ have a common fixed point.
\end{enumerate}
\end{lemma}

\begin{proof}
The first point is Corollary 2.5 in \cite{PECAbstract}.  For the second point, note that $N.U=\langle N,U \rangle$ is boundedly generated by $N$ and $U$. Then a fixed point $x_0$ for $N.U$ is fixed by both $N$ and $U$.
\end{proof}

\subsection{Kac--Moody groups} Let $n \in \NN$ and $A=(a_{ij}) \in \mathbf{Z}^{n \times n}$ be a generalized Cartan matrix, i.e. $a_{ii}=2$ and $a_{ij} \leq 0$ for $i \neq j$. Let $\Lambda$ be a free $\ZZ$-module of finite rank and $c_i \in \Lambda$, $h_i \in \Lambda^\vee=\Hom(\Lambda,\ZZ)$ for $i=1, \ldots, n$ be such that $h_i(c_j)=a_{ij}$. Then the quadruple $\mathcal D=(A,\Lambda, (c_i),(h_i))$ is called a \textbf{Kac--Moody root datum}. If $\Lambda$ is spanned by the $c_i$, the datum $\mathcal D$ is called \textbf{adjoint}, while if $\Lambda^\vee$ is generated by the $h_i$ it is called \textbf{simply connected}.\\

Let $\mathcal G_D$ denote the \textbf{constructive Tits functor} associated to $\mathcal D$ \cite{Tits87}. This is a functor from the category of commutative rings with 1 to the category of groups satisfying certain natural properties, the most important of which is the existence of the adjoint representation over the complex numbers. This is to say that there is a homomorphism $\Ad\colon \mathcal G_{\mathcal D}(\CC) \to \Aut(\mathfrak g_A)$, where $\mathfrak g_A$ denotes the complex Kac--Moody algebra associated to $A$, whose kernel is central in $\mathcal G_{\mathcal D}(\CC)$. For a field $k$, the group $\mathcal G_\mathcal D(k)$ is called a
(split) \textbf{Kac--Moody group} over $k$. It can be given explicitly by generators and relations, see e.g. \cite{CarboneErshovRitter}.\\
\begin{example} 
\begin{enumerate}
\item Let $\mathcal D$ be the root datum associated to a reductive algebraic group $G$ over an algebraically closed field $k$. Then 
$\mathcal G_\mathcal D(k)$ coincides with $G(k)$. In this case, $\mathcal G_\mathcal D$ is said to be \textbf{of finite type}. 
\item Let $A$ be a generalized Cartan matrix of affine type and $\mathcal D$ the simply-connected root datum associated to $A$. Then $\mathcal G_\DD(k)$ is called an \textbf{affine Kac--Moody group} and can be identified with a central extension of a loop group. For example, the group $\SL_n(k[t,t^{-1}])$
is an affine Kac--Moody group.
\end{enumerate}
\end{example}

Let $k$ be a field and $G:=\mathcal G_\mathcal D(k)$ a Kac--Moody group over $k$. Then there are subgroups $B_+,B_-\leq G$ such that $(G,B_+,B_-,S)$ is a twin BN-pair for $G$ (cf. \cite[Chapter 8]{AbramenkoBrown}). Let $(W,S)$ denote the associated Coxeter group.

There are Bruhat (resp. Birkhoff) decompositions
$$G= \bigcup_{w \in W} B_+ w B_+=\bigcup_{w \in W} B_- w B_-=\bigcup_{w \in W} B_+wB_-=\bigcup_{w \in W} B_-wB_+.$$

For $\epsilon \in \{\pm 1\}$, the sets $G/B_{\epsilon}$ respectively form the set of chambers of a building $\Delta_\epsilon$ of type $(W,S)$. This simplicial complex has a geometric realization $|\Delta_\epsilon|$ which makes it into a complete CAT(0)-space (cf. \cite[Chapter 12]{AbramenkoBrown}). The group $G$ acts on $|\Delta_\epsilon|$ by isometries. The buildings $\Delta_\epsilon$ are tied together by a codistance $\delta^*$, which makes $\Delta=(\Delta_+,\Delta_-)$ into a twin building (cf. \cite[Chapter 8]{AbramenkoBrown}).\\ 
A subset $\Omega \subseteq \Delta$ which is contained in a twin apartment $\mA=(\mA_+,\mA_-)$ is called {\bf balanced} if 
$\Omega \cap \mA_+ \neq \emptyset \neq \Omega \cap \mA_-$ and $\Omega$ is contained in the union of a finite number of spherical residues.
Note that if $x_\pm \in |\Delta_\pm|$, then $\Omega=\{x_+,x_-\}$ is a balanced set. Indeed, if $C_\pm$ are chambers containing $x_\pm$ respectively, then there is 
a twin apartment $\mA$ containing $C_\pm$ (hence $x_\pm$), and the points $x_\pm$ correspond to spherical residues. 

A subgroup $P_+$ (resp. $P_-$) which contains a conjugate of $B_+$ (resp. $B_-$) is called a positive (resp. negative) \textbf{parabolic subgroup}.
In this case there is a subset $J\subseteq S$ such that $P_+=B_+W_JB_+=:P_{+,J}$, similarly for $P_-$. A positive parabolic $P_+$ and a negative parabolic $P_-$ are called \textbf{opposite} if there is $g \in G$ and $J \subseteq S$ such that $P_\pm=gP_{J,\pm}g^{-1}.$ The intersection $P_+ \cap P_-$ of two opposite parabolics is called a \textbf{Levi factor}. It can be explicitly described as folllows: If $(T,(U_\alpha)_{\alpha \in \Phi(W,S)})$ denotes the standard twin root datum for $G$, then $L_J:=P_{J,+} \cap P_{J,-}=T\cdot \langle U_\alpha: \alpha \in \Phi(W_J,J) \rangle$. The twin building associated to $L_J$ embeds naturally into $\Delta$.\\
For $\epsilon \in \{\pm 1\}$, a parabolic subgroup $P_{\epsilon}$ is called \textbf{of finite type} if it is contained in a finite number of double $B_\epsilon$-cosets. If $P=gP_{J,\epsilon} g^{-1}$, $P$ is of finite type if and only if $W_J$ is finite. \\
Two points $x_\pm \in |\Delta_\pm|$ are called \textbf{opposite} if $\Stab_G x_\pm$ are opposite parabolic subgroups of finite type.

\subsection{Bounded subgroups of Kac--Moody groups}
Let $G:=\mathcal G_\mathcal D(k)$ be a Kac--Moody group.
A subgroup $H \leq G$ is called \textbf{bounded} if $H \leq P_+ \cap P_-$ for a positive parabolic $P_+$ and a negative parabolic $P_-$, both of which are of finite type.\\
An important characterization is as follows: $H$ is a bounded subgroup of $G$ if and only if its actions on both $|\Delta_+|$ and on $|\Delta_-|$ are bounded \cite[Theorem 10.2.2]{Remy}.\\

Via the adjoint representation, it is possible to endow bounded subgroups with the structure of an algebraic group. For the following construction due to R\'emy, we refer to \cite[Section 10.3.1]{Remy}, the notation of which we use freely. \\
Let $\Omega \subseteq \Delta$ be a balanced subset which is contained in the standard twin apartment $\mA$ of $\Delta$. 
Let $\bar{k}$ be an algebraic closure of $k$. 
Let $\mcL_\DD:=\lieg_\DD \cap \UU_\DD$, where $\UU_\DD$ is the $\ZZ$-form of the universal enveloping algebra of the associated Kac--Moody algebra. Then $\mcL$ has a grading: $\mcL_\DD=\mcL_0 \oplus \bigoplus_{a \in \Phi}\mcL_a.$

Let $\Delta(\Omega):=\{a \in \Phi: \Omega \subseteq a\}$, $\Delta^u(\Omega):=\{a \in \Phi: \Omega \subseteq a, \Omega \subsetneq \partial a\}$ and let $\Delta^m(\Omega):=\{a \in \Phi: \Omega \subseteq \partial a\}.$
Here the roots are viewed as half-spaces.
Write $L:=T\langle U_\alpha: \alpha \in \Delta^m(\Omega)\rangle $ and $U:=\langle U_\alpha: \alpha \in \Delta^u(\Omega)\rangle.$ 
 
\begin{proposition} \label{boundedsubspace}
Let $W=W_\Omega$ be the smallest $Q$-graded subspace of $(\mcL_\DD)_{\bar{K}}$ with the following properties:
\begin{enumerate}
\item $W$ contains $(\mcL_0)_{\bar{K}}$ and $(\mcL_a)_{\bar{K}}$ for all $a \in \Delta(\Omega)$.
\item The $Q$-support of $W$ contains $-\Delta^u(\Omega).$
\item $W$ is stable under $H:=\Fix \Omega.$
\end{enumerate}
Then the following properties hold:
\begin{enumerate}
\item 
$W$ is finite-dimensional and the kernel of $\Ad\colon H \to \Ad H|_W$ is precisely the center of $H$. 
\item Let $\bar{H}$ (resp. $\bar{T}, \bar{L}, \bar{U}$) denote the Zariski-closure of $\Ad H|_W$ (resp. $\Ad T|_W$, $\Ad L|_W$, $\Ad U|_W$).
Then $\bar{L}$ is a connected reductive $K$-group, $\bar{T}$ is a maximal torus of $\bar{L}$, $\bar{U}$ is unipotent and $\bar{H}=\bar{L} \ltimes \bar{U}$ is a Levi decomposition. \index{Levi decomposition}
\end{enumerate}
\end{proposition}

\begin{proof}
This is \cite[Lemma 10.3.1, Proposition 10.3.6]{Remy}.
\end{proof}

An important fact about Kac-Moody groups is that checking if a subgroup is bounded can be done by using functoriality. For the theory of almost split Kac--Moody groups see \cite{Remy}. 
\begin{lemma} \label{subbuildingbounded} Let $E$ be a field and let $k$ be a subfield of $E.$ Let $\GG_\DD$ be a Tits functor and let $G:=\GG_\DD(E).$
Let $H$ be one of the following subgroups of $G$:

\begin{enumerate}
\item $H=G(k)$, where $G(k)$ is an almost split Kac--Moody $k$-group obtained from $\GG_\DD(E)$ by Galois descent.
\item $H=L(k)$, where $L$ is a Levi factor of $G(k)$.
\end{enumerate}
Let $U\leq H$. Then $U$ is bounded with respect to $H$ if and only if $U$ is bounded with respect to $G$.
\end{lemma}

\begin{proof}
The twin building $\Delta(H)$ associated to $H$ embeds as a closed convex subbuilding into the twin building $\Delta(G)$ associated to $G$. If $U$ fixes points in both halves of $\Delta(H)$, the same points serve as fixed points in $\Delta(G)$. If $U\leq H$ fixes a point in $\Delta(G)$, it must also fix a point in $\Delta(H)$ since it leaves this closed convex set invariant, cf. \cite[Proposition II.6.2 (4)]{BridsonHaefliger}. 
\end{proof}

The following two propositions collect standard facts of the theory of algebraic groups, which have a geometric meaning in the context of twin buildings. 
\begin{proposition} \label{Hochschild}
Let  $k$ be a field of characteristic 0 and let $G$ be an algebraic group over $k$. Let $G_u$ denote the unipotent radical of $G$. Let $P$ be a linearly reductive subgroup of $G$ such that $G=G_uP$ and let $Q$ be any linearly reductive subgroup of $G$. Then there exists $t \in G_u(k)$ such that $tQ(k)t^{-1} \subseteq P(k).$
\end{proposition}

\begin{proof}
This is Proposition VIII.4.2 in \cite{MR620024}.  
\end{proof}
\begin{proposition}\label{SL2Rational}
Let $k$ be a field of characteristic 0 and $\overline{k}$ an algebraic closure of $k$. Let $n \in \NN$ and let $\varphi\colon \SL_2(\QQ) \to \GL_n(k)$ be a group homomorphism. Then there is a homomorphism of algebraic groups $\psi\colon \SL_2 \to \GL_n$ defined over $k$ such that $\psi|_{\SL_2(\QQ)}=\varphi.$  
\end{proposition}
\begin{proof} 
This is in \cite[p. 343]{MR822247}. Another proof is given in \cite[Lemma 5.9]{PECAbstract}.
\end{proof}

The following proposition is a refinement of an important result by Caprace.
\begin{proposition}\label{sl2oppoints}
Let $k$ be a field and let $G$ be an almost split Kac-Moody group over $k$. Let $\varphi\colon \SL_2(\QQ) \to G(k)$ be a homomorphism. 
Then $\im \varphi$ is bounded and fixes two opposite points of $\Delta(G).$
\end{proposition}

\begin{proof}
Let $E$ be a field over which $G$ splits and let $\iota\colon G(k) \to G(E)$ denote the canonical inclusion.
By \cite[Corollary 5.8]{PECAbstract}, $\iota \circ \varphi$ has bounded image and if $\varphi$ is non-trivial (which we may assume), $E$ has characteristic 0. By Lemma \ref{subbuildingbounded}, $\im \varphi$ fixes points $x,y$ in both halves of $\Delta(G(k)).$ Then $\Omega:=\{x,y\}$ is balanced. Let $\psi:=\Ad_\Omega \circ \varphi\colon \SL_2(\QQ) \to \GL(W_\Omega).$ Then $\psi$, although a priori only an abstract representation, is in fact rational by the previous proposition. It follows that the Zariski closure $C$ of $\psi(\SL_2(\QQ))$ is a reductive group. By Proposition \ref{Hochschild}, $C$ can be conjugated inside a Levi factor of $\Ad_\Omega(\Fix \Omega)$ by an element of $G(k)$. This Levi factor is precisely the stabiliser of two opposite points, from which the claim follows. 
\end{proof}

\section{Unbounded algebraic subgroups}
Let $k$ be a field with infinite transcendence degree over its prime field and let $G_1$ be a connected $k$-isotropic algebraic $k$-group. In this section, we will give a construction of a homomorphism $\varphi\colon G_1(k) \to G_2(k)$, where $G_2(k)$ is a certain Kac--Moody group over $k$, with unbounded image.  \\

In the case of affine Kac--Moody groups, there is an easy criterion to check if a subgroup is bounded. Let $k$ be a field and let $\deg_t$, $\deg_{t^{-1}}$ denote the valuations of $k(t)$ with $t$ and $t^{-1}$ as uniformizing parameter, respectively.
For an element 
$$f=\sum_{i=N_0}^{N_1} a_it^{i} \in k[t,t^{-1}]$$
with $a_{N_0} \neq 0, a_{N_1} \neq 0$ this means that $\deg_t(f)=N_0$.
Since $\deg_{t^{-1}}(f)=\deg_t f(\frac 1 t)$ it follows that $\deg_{t^{-1}}(f)=-N_1$.\\
For a matrix $g=(g_{ij}) \in \SL_n(k[t,t^{-1}])$, let $\deg_t(g):=\displaystyle \max_{i,j} \deg_t g_{ij}$ and similarly let $\deg_{t^{-1}}(g):=\displaystyle \max_{i,j} \deg_{t^{-1}} g_{ij}$. \\
With this notation, we have the following:
\begin{proposition}\label{boundedinaffine}
A subgroup $U$ of $\SL_n(k[t,t^{-1}])$ is bounded if and only if there exists $N \in \NN$ such that $|\deg_t(g)|\leq N$ and $|\deg_{t^{-1}}(g)|\leq N$ for all $g \in U.$
\end{proposition}
\begin{proof}
The positive (resp. negative) half of the twin building can be identified with the affine building of $\SL_n(k(t))$ with the discrete valuation $v$ (resp. $v'$) having $t$ (resp. $\frac{1}{t}$) as a uniformizing parameter. By \cite[Ex. 11.40] {AbramenkoBrown} $U$ is bounded if and only if there is an upper bound on the absolute values of the matrix entries, which amounts to the claim. 
\end{proof}

Let $k$ be a field and $G_1$ a $k$-isotropic reductive algebraic $k$-group.
Let $\mathfrak g_1$ denote the Lie algebra of $G_1$. Then $G_1(k)$ acts on $\mathfrak g_1(k)$, which is as an abelian group isomorphic to $k^n$, via the adjoint representation. 
Let $G_2:=\SL_{n+1}(k[t,t^{-1}])$ and consider the following subgroup: 
$$V=\left\{
    \begin{pmatrix}
     &       & \vdots \\
     & \Ad g & v \\
     &       & \vdots    \\
     &       & 1   \\
    \end{pmatrix}: g \in G_1(k), v \in (k[t,t^{-1}])^n \right\}.$$

Note that $\det \Ad g=1$ is automatic since $G_1$ is reductive. Indeed, it suffices to check this over an algebraically closed field. Then $G_1$ is the almost direct product of its center $C$ and its derived group $[G_1,G_1]$. Now $\Ad C$ is trivial, whereas $\det \Ad g=1$ for each element $g$ of the derived group.

\begin{lemma}\label{structureofgroup}
As an abstract group, $V$ is isomorphic to 
$$ \Ad G_1(k) \ltimes \bigoplus_{i \in \ZZ} \mathfrak g_1(k),$$
where the action of $\Ad G_1$ on each summand of the direct sum is the natural one.
\end{lemma}
\begin{proof}
 This is clear since $(k[t,t^{-1}],+)^n \cong \bigoplus_{i \in \ZZ} k^n$ by decomposing an element into its homogeneous components, and this decomposition is preserved under the action of $\Ad g$. The given matrix representation is just the standard one for a semidirect product of two linear groups.
\end{proof}

Until the end of this section assume that the transcendence degree of $k$ over its prime field $F$ is infinite. Let $T={(t_i)_{i \in \ZZ}}$ be an infinite set of algebraically independent elements. Complete $T$ by a set $T'$ to a transcendence basis for $k$ over $F$. For each $i \in \ZZ$, consider the derivation $\delta_i\colon k \to k$ obtained 
by extending the zero-derivation on $F(T'\cup T \backslash \{t_i\})$ to $k$ by setting $\delta_i(t_i)=1$, i.e. $\delta_i$ can be thought of as a partial derivative with respect to $t_i.$ By the basic theory of derivations (cf. \cite{Jacobson}), this is clearly possible. \\ 
Since for each $x \in k$ there is a finite set $I_x$ such that $x$ is contained in an algebraic extension of $F(t_i: i \in I_x)$, it follows that 
for each $x \in k$, $\delta_i(x)\neq 0$ for only finitely many $i.$ \\

The following observation is due to Borel--Tits \cite[Example 8.18 b)]{MR0316587} and is elaborated on in Lifschitz--Rapinchuk \cite{MR1844711}. 
For a closed subgroup $G$ of $\GL_n(k)$ and a derivation $\delta\colon k \to k$, the matrix $g^{-1}\cdot \delta(g)$ is an element of $\lieg=\Lie G$, where $\delta(g)$ is the matrix obtained by applying $\delta$ to all entries of $g$. Furthermore, the mapping $\varphi_\delta\colon G \to G \ltimes \lieg, g \mapsto (g, g^{-1}\cdot \delta g)$ is a group homomorphism. Here $G$ acts on $\lieg$ via the adjoint representation.\\

This construction can be used to describe homomorphisms with unbounded image.

\begin{theorem}\label{infinitetranscendenceunbounded}
Let $k$ be a field with infinite transcendence degree over its prime field. Let $G_1$ be a $k$-isotropic reductive $k$-group and let $n:=\dim G_1$. Let $G_2$ be a Kac-Moody group whose diagram has a subdiagram of type $\widetilde{A_n}$ such that the derived group of the corresponding Levi factor $H(k)$ is isomorphic to $\SL_{n+1}(k[t,t^{-1}]).$\\
Let $V=\Ad G_1(k) \ltimes \bigoplus_{i \in \ZZ} \mathfrak g_1(k) \leq H(k)$ be as above and let $(\delta_i)_{i \in \ZZ}$ be derivations as above.
Then the mapping
$$\varphi\colon G_1(k) \to V\leq G_2, g \mapsto (\Ad g, (g^{-1}\cdot \delta_i(g))_{i \in \ZZ})$$
is a group homomorphism with unbounded image.  
\end{theorem}

\begin{proof}
Due to the construction of the derivations it follows that for each $g \in G_1(k)$, $(g^{-1}\cdot \delta_i(g))\neq 0$ for only finitely many $i \in \ZZ$, hence $\varphi$ is well-defined. Since each $\varphi_{\delta_i}$ is a homomorphism, so is $\varphi.$\\
Let $T(k)\leq G_1(k)$ be a $k$-split torus. Then for each $i$ there is an element $g_i \in T(k)\cong (k^\times)^r$ such that $g^{-1}\cdot \delta_i(g)\neq 0$. 
This translates via Proposition \ref{structureofgroup} to the fact that $\varphi(g_i)$ has a matrix entry which has a homogeneous component of degree $i$. In particular, the degrees of the matrix entries of $\varphi(G)$ are unbounded, which proves via Proposition \ref{boundedinaffine} that $\varphi$ has unbounded image in $H(k).$ By Lemma \ref{subbuildingbounded}, $\im \varphi$ is unbounded in $G$, too. 
\end{proof}

\begin{remark} If a field $k$ is uncountable, it is of infinite transcendence degree over its prime field. In particular, a local field has infinite transcendence degree over its prime field, as follows from the classification of local fields (cf. \cite{MR554237}). 
\end{remark}

This result is interesting in a different context, too. Following Farb \cite{Farb}, a group $G$ is said to have property FA$_n$ if every $G$-action by cellular isometries on an $n$-dimensional CAT(0)-complex has a global fixed point. \\
Note that the geometric realization of $\SL_{n+1}(k[t,t^{-1}])$ has dimension $n$, since the apartments are tessellations of $\RR^n$. 
This implies the following corollary.

\begin{corollary}
Let $k$ be a field with infinite transcendence degree over its prime field and let $G$ be a reductive isotropic $k$-group of dimension $n$.
Then $G(k)$ does not have FA$_n$. 
\end{corollary}

The following discussion, leading up to Corollary \ref{cantcheck}, aims to make precise the informal statement that the space of quasi-morphisms of a Kac--Moody group $G$, although infinite-dimensional, cannot be used to check whether a subgroup is bounded or not. \\

For a group $G$, a map $\varphi\colon G \to \RR$ is called a {\bf quasi-morphism} if it satisfies
$$\sup_{g,h \in G} |\varphi(gh)-\varphi(g)-\varphi(h)|<\infty.$$
Let QH$(G)$ denote the real vectorspace of all quasi-morphisms of $G$.
Then $l^\infty(G)$, the space of bounded real-valued functions on $G$, and $\Hom(G,\RR)$ are subspaces of QH$(G)$. 
Let $\widetilde{\QH}(G):=\QH(G)/(l^\infty(G) \oplus \Hom(G,\RR))$ denote the space of non-trivial quasi-morphisms of $G$. 

\begin{theorem}\label{classificationquasihom}
Let $k$ be a field and let $G:=\GG_\DD(k)$ be a split Kac-Moody group such that the Weyl group $(W,S)$ of $G$ is irreducible and neither spherical nor affine. Then $\widetilde{\QH}(G)$ is infinite-dimensional. 
\end{theorem}

\begin{proof}
This is \cite[Theorem 1.1]{CapraceFujiwara}.
\end{proof}

Let $k$ be a local field and let $G$ be a connected simply connected almost simple algebraic group defined over $k$ of $k$-rank $\geq 2.$
By \cite[Lemma 6.1]{BurgerMonod1999}, any {\it continuous} quasi-morphism $f\colon G(k) \to \RR$ is trivial. When one restricts to $G=\SL_n(k)$ for $n \geq 3$, it is possible to drop the continuity assumption.
\begin{proposition} \label{sl3trivial}
Let $k$ be a field and let $n\geq 3.$ Then any quasi-morphism $f\colon \SL_n(k) \to \RR$ is bounded.
\end{proposition}
\begin{proof}
The group $\SL_n(k)$ is boundedly generated by its root groups $U_\alpha$, so it suffices to show that $f$ is bounded on each $U_\alpha.$
Since $n\geq 3$, any two elements $a, b \in U_\alpha \backslash \{1\}$ are conjugate via some diagonal matrix. 
Since $f$ is bounded on conjugacy classes, the claim follows.
\end{proof}

For a group homomorphism $\varphi: H \to G$ there is a pull-back
$\varphi^*\colon \QH(G) \to \QH(H)$ given by $\varphi^*(f)(x):=f(\varphi(x))$. 

\begin{corollary}\label{cantcheck}
Let $k$ be a field with infinite transcendence degree over its prime field. Let $G:=\GG_\DD(k)$ be a split Kac-Moody group such that $G$ contains $\SL_4(k[t,t^{-1}])$ as the derived group of a Levi factor and such that the Weyl group of $G$ is irreducible and not affine. Then $\QH(G)$ is infinite-dimensional.\\
There is a homomorphism
$$\varphi\colon \SL_3(k) \to G$$
such that $\im \varphi$ is unbounded yet $|f(\im \varphi)|<K_f$ for each quasi-morphism $f \in \QH(G)$ and a constant $K_f$ depending on $f$.
\end{corollary}

\begin{proof} Note first that such a Kac-Moody group $G$ exists, as it suffices to extend the root datum associated to $\SL_4(k[t,t^{-1}])$ to make the associated Cartan matrix to be of indefinite type. Then the first statement follows from Theorem \ref{classificationquasihom}.\\
Let $\varphi\colon \SL_3(k) \to G$ be the homomorphism with unbounded image constructed in Theorem \ref{infinitetranscendenceunbounded}. For $f \in \QH(G)$, the pull-back $\varphi^* f$ is bounded by Proposition \ref{sl3trivial}, from which the claim follows.
\end{proof}

\section{The case of number fields}
In this section we show that embeddings of split Chevalley groups over number fields into Kac--Moody groups always have bounded image.

\begin{lemma}
Let $L|\QQ$ be a finite field extension and let $x \in L$ be a primitive element, i.e. $L=\QQ(x).$
Then for $n \in \NN$, there exists  $y \in L$ such that $y,y^2, \ldots, y^n$ all are primitive elements. 
\end{lemma}

\begin{proof}
Set $k:=n!$ and observe that it is enough to find an element $y$ such that $y^k$ is primitive, since then for any $p$ dividing $k$ we clearly have $\QQ(y^p) \supseteq \QQ(y^k)=L.$
Consider the sequence $y_i:=x+i,  i \in \NN$ and set $z_i:=y_i{}^k.$ Since there are only finitely many intermediate fields between $L$ and $\QQ$, by the pigeonhole principle there is a subsequence $(z_{i_l})$ such that $\QQ(z_{i_r})=\QQ(z_{i_s})$ for all $i_r, i_s$. \\
 Note that in the polynomial ring $\QQ[t]$ for pairwise different $a_1, \ldots, a_{k+1} \in \NN$ the polynomials $(t-a_i)^k$ form a basis of the subspace of polynomials of degree $\leq k$, since the coefficient vectors form a Vandermonde matrix (up to scaling). In particular, $t$ lies in their $\QQ$-span. Applying this analysis to $\QQ(z_{i_r})$ then shows that this field contains $x$ and hence equals $L.$ 
\end{proof}

\begin{lemma}
Let $L$ be a number field of degree $n$. Let $x \in L$ be such that $x^2$ is again a primitive element. 
Then for $d:=\diag(x,x^{-1})$, $\SL_2(L)$ is boundedly generated by $d$ and $\SL_2(\QQ).$  
\end{lemma}
\begin{proof}
Set $u_+(r):=\left( \begin{smallmatrix}
              1 & r \\
              0 & 1
             \end{smallmatrix} \right)$ 
and note that $d^{i} u_+(r) d^{-i}=u_+(rx^{2i}).$ This implies that $U_L:=\{u_+(l): l \in L\}$ is generated by the subgroups $d^{i}U_{\QQ}d^{-i}, i=0\ldots,n-1.$ 

Now $\SL_2(L)$ is boundedly generated by $U_L$ and $sU_Ls^{-1}$, where $s:=\left( \begin{smallmatrix}
              0 & -1 \\
              1 & 0
             \end{smallmatrix} \right)$ 
, i.e. by $2n$ conjugates of $U_{\QQ}$.
\end{proof}

\begin{corollary} \label{numberfieldbounded} Let $\GG_\DD$ be a constructive Tits functor and $K$ a field. Then for any number field $L$, every homomorphism $\SL_2(L) \to \GG_\DD(K)$ has bounded image.
\end{corollary}
\begin{proof}
Since $U_\QQ \subseteq \SL_2(\QQ)$ has fixed points in both halves of the twin building by Proposition \ref{sl2oppoints}, so does every conjugate of $U_\QQ$. By the preceding lemma and Lemma \ref{boundedgenerationlemma}, the claim follows. 
\end{proof}
\begin{remark} 
\begin{enumerate}
\item Using exactly the same reasoning, one can show: If $L|K$ is a finite separable extension, char $K\neq 2$ and $\varphi\colon \SL_2(L) \to \GG_\DD(F)$ is a homomorphism whose restriction to $\SL_2(K)$ has bounded image, then $\im \varphi$ also is bounded. 
\item As already remarked in \cite{PECAbstract}, via bounded generation the above implies that any homomorphism of a split Chevalley group over a number field has bounded image.
\end{enumerate}
\end{remark}

The fixed point set of $\SL_2(L)$ might be smaller than that of its subgroup $\SL_2(\QQ)$, as the next example shows.
\begin{example} Let $L$ be a number field of degree $d>1$ and let $\sigma_1, \sigma_2\colon L \to \CC$
be two different embeddings. Let $G=\SL_2$ and consider the homomorphism
$$\varphi\colon G(L) \to G(\CC) \times G(\CC), g \mapsto (\sigma_1(g), \sigma_2(g)).$$ 
Postcomposing with the standard inclusion $G(\CC) \times G(\CC) \to \SL_4(\CC)$, $G(L)$ then acts on the spherical building associated to $\SL_4(\CC)$. While $G(\QQ)$ fixes the residue $R$ corresponding to the subspace $U:=\langle (1,0,1,0), (0,1,0,1)\rangle$, the group $G(L)$ does not fix it. \\
In particular, $\varphi(G(\QQ))$ is not Zariski dense in the closure of $\varphi(G(L))$: the closure of $\varphi(G(\QQ))$ is the diagonal subgroup
of the closure of $\varphi(G(L))$, which is $\SL_2(\CC) \times \SL_2(\CC)$. In particular, the dimension of the Zariski closure increases.\\
This observation leads to the proof of Theorem \ref{infinitealgebraicbounded} in the next section. 
\end{example}

\section{The case of infinite algebraic extensions}
In the following we address the question of infinite algebraic extensions. 
We use Margulis's rigidity result that any abstract representation of $\SL_2(R)$, where $R$ is a certain subring of a number field, has semisimple Zariski closure. \\

\index{lattice}
\index{lattice!irreducible}
Recall first that for a locally compact topological group $G$, a {\bf lattice} $\Gamma$ is a discrete subgroup of $G$ such that $G/\Gamma$ has a finite invariant measure. For locally compact topological groups $G_1, \ldots, G_n$, a lattice $\Gamma \leq G_1 \times \ldots \times G_n$ is said to be {\bf irreducible} if the projection of $\Gamma$ on each factor $G_i$ is dense in $G_i.$ By \cite[Corollary 5.21]{Raghunathan} this is equivalent to the condition that no subgroup of finite index is a direct product of two non-trivial normal subgroups.\\

The reference for the following paragraph is \cite[Introduction]{Margulis}.\\
Let $L$ be a number field and let $R$ be the set of all (inequivalent) valuations of $L$. Let $R_\infty$ denote the set of archimedean valuations of $L$, and for each $v\in R$ let $L_v$ be the completion of $L$ with respect to $v.$ Let $|\cdot|_v$ denote the absolute value associated with $v$.\\
Let $S \subseteq R$ be a finite subset containing all archimedean valuations and suppose $|S|\geq 2.$ Let $L(S):=\{x \in L: |x|_v \leq 1$ for all non-archimedean $v \in R \backslash S\}$ be the ring of $S$-integral elements of $L$.

\begin{theorem}[Borel--Harish-Chandra]
Let $G=\SL_2$. Then $\Gamma:=G(L(S))$ is an irreducible lattice in $G_S:=\prod_{v \in S} G(L_v)$.
\end{theorem}

\begin{proof} 
This is the special case $G=\SL_2$ of the result quoted in \cite[Page 1]{Margulis}.
\end{proof}
The following lemma combines several standard facts about lattices. 
\begin{lemma} \label{Gamma0}
Let $\Gamma=G(L(S))$ be as above. Then there is a subgroup $\Gamma_0 \leq \Gamma$ of finite index such that $\Gamma_0$ has trivial center and still is an irreducible lattice in $G_S.$ 
\end{lemma}

\begin{proof}
By Selberg's lemma for finitely generated linear groups it follows that $\Gamma$ is virtually torsion free, i.e. there is a subgroup $\Gamma_0\leq \Gamma$ of finite index which is torsion free.\\
Since $\Gamma_0\leq \Gamma$ and $\vol (G/\Gamma_0)=[\Gamma : \Gamma_0]\vol(G/\Gamma)$, $\Gamma_0$ is a lattice. \\
By the second condition for irreducibility given above, it is clear that $\Gamma_0$ is irreducible, as a direct product of finite index in $\Gamma_0$ would be of finite index in $\Gamma$.
By Margulis's Normal Subgroup Theorem, any central element in $\Gamma_0$ has finite order, hence must be trivial.
\end{proof}

In this setting, Margulis's Rigidity Theorem takes the following form.
\begin{theorem}\label{MargulisThm} Let $G_S=\prod_{v \in S}\SL_2(L_v)$ be as above and let $\Gamma \leq G_S$ be an irreducible lattice. 
Let $K$ be a field of characteristic 0, $F$ an algebraic group defined over $K$ and $\varphi\colon \Gamma \to F(K)$ a homomorphism. Then the Zariski closure of $\varphi(\Gamma)$ is a semisimple algebraic group defined over $K$.
\end{theorem}
\begin{proof}
This is Theorem 3 of \cite[Introduction]{Margulis} adapted to the present situation. 
\end{proof}
\begin{remark} The Zariski closure of $\varphi(\Gamma)$ need not be connected since in general $\Gamma$ will admit finite quotients.
\end{remark}

To prove Theorem \ref{infinitealgebraicbounded}, we will apply Theorem \ref{MargulisThm} to the lattice $\Gamma_0$ provided by Lemma \ref{Gamma0}. Before proceeding to the proof of this theorem, we need some lemmas on bounded subgroups. \\

Let $U\leq \GG(k)$ be a bounded subgroup of a split Kac-Moody group $\GG(k).$ Let $\Delta^U$ denote the fixed point set of the $U$-action on the CAT(0)-realization of the twin building. For $x \in \Delta_+^U$ and $y \in \Delta_-^U$ let $\Omega(x,y):=\{x,y\}$ and $\mA$ a twin apartment containing $\Omega$. 

Let $W:=W_{\Omega(x,y),\mA}$ of $\UU_\DD$ be the finite-dimensional $\Fix(\Omega(x,y))$-invariant subspace provided by Proposition \ref{boundedsubspace}.
Let $U_{x,y,\mA}$ denote the Zariski closure of $\Ad U|_{W_{\Omega(x,y),\mA}}$. 

\begin{remark}
The group $U_{x,y,\mA}$ depends in general on the choice of $x$ and $y$. For instance, let $U=T$ be the standard torus which fixes the standard twin apartment $\mA=(\mA_+,\mA_-).$ Let $C_+$ be the standard positive chamber and $C_-$ the standard negative chamber. \\
Choose $x_0$ in the interior of $C_+$ and $y_0$ in the interior of $C_-.$ Then $\Delta({x,y})=\emptyset$ and by the calculation of $\ker \Ad$ it follows that $U_{x_0,y_0,\mA}=1.$\\
On the other hand, let $\alpha$ be a simple root and let $x_1 \in \mA_+$ be in the interior of $\partial \alpha$ but not in the interior of any other $\partial \beta$, where $\beta$ is a positive simple root distinct from $\alpha$. Let $y \in \mA_-$ be chosen similarly. Then 
$\Delta({x_1,y_1})=\Delta^m({x_1,y_1})=\{\alpha, -\alpha\}.$
By the calculation of $\ker \Ad$ it follows that $U_{x_1,y_1,\mA}$ is a 1-dimensional torus.
\end{remark}

The point of this remark is that one has to restrict to bounded subgroups which have trivial center. 

\begin{theorem} Let $k$ be a field and let $G:=\GG_\DD(k)$ be a split Kac-Moody group. Let $U\leq G$ be a bounded subgroup with trivial center and let $x,x' \in \Delta_+, y,y' \in \Delta_-$ be $U$-fixed points. Let $\mA$ be a twin apartment which contains $x,y$ and $\mA'$ a twin apartment which contains $x',y'.$ 
Let $U_{x,y,\mA}$ denote the Zariski closure of $\Ad U|_{W_{\Omega(x,y),\mA}}$ and similarly $U_{x',y',\mA'}$ for $x',y',\mA'.$ \\
Then $U_{x,y,\mA}$ and $U_{x',y',\mA'}$ are isomorphic as algebraic groups.
\end{theorem}

\begin{proof}
Note first that if $x_0,y_0$ are two points fixed by $U$ and $\mA_0$ is a twin apartment containing them, then the restriction of $\Ad U$ to $W_{x_0,y_0,\mA_0}$ is injective: the kernel of this restriction is central in $\Fix_G\{x_0,y_0\}$, so in particular its intersection with $U$ must be central in $U$, but $U$ has trivial center.\\

{\it Step 1. If $\mA''$ is another apartment containing $x$ and $y$, then $U_{x,y,\mA} \cong U_{x,y,\mA''}$.}\\
The group $\Fix_G \{x,y\}$ acts transitively on the set of apartments containing $x$ and $y$ by \cite[Proposition 10.4.4 (iii)]{Remy}. Let $g \in \Fix_G\{x,y\}$ be such that $g\mA=\mA''.$ Then
$\Ad g$ conjugates $W_{x,y,\mA}$ to $W_{x,y,\mA''}$ and hence $U_{x,y,\mA}$ onto $U_{x,y,\mA''}$.\\

Since $U_{x,y,\mA}$ is independent from the twin apartment $\mA$ containing $x,y$, we can write $U_{x,y}:=U_{x,y,\mA}$ unambiguously.\\

{\it Step 2. It suffices to prove that $U_{x,y} \cong U_{x',y}$, where
$x,x',y$ are contained in a common apartment and $x,x'$ are contained in a common chamber}\\
Let $\rho\colon [0,r] \to \Delta_+, \rho(0)=x,\rho(r)=x'$ be a geodesic, where $r:=d(x,x').$ Then $\im \rho$ is fixed by $H$.\\
Let $0=i_0<i_1< \ldots<i_k=r$ be such that $\rho([i_j,i_{j+1}])$ is contained in a chamber $C_j$, $j=0, \ldots, k-1$. Let $x_j:=\rho(i_j).$
Suppose it is already proven that $U_{x_j,y} \cong U_{x_{j+1},y}.$
Since $x_0=x, x_k=x'$ it then follows that $U_{x,y} \cong U_{x',y}$,
and arguing similarly for $y$ it follows that $U_{x',y} \cong U_{x',y'}.$\\

{\it Step 3. Conclusion.}~Let $F_I,F_{I'}$ be facets containing $x,x'$ and maximal with this property. Then $H$ fixes both $F_I$ and $F_I'$
and hence $F_{I \cap I'}.$ Replace the geodesic from $x$ to $x'$ by the union of a geodesic from $x$ to a point $z$ in $F_{I \cap I'}$ and a geodesic from $z$ to $x'$. This allows to assume that $I \subseteq I'$ or $I' \subseteq I$, without loss of generality assume that $I \subseteq I'.$ \\
With the notations from \cite[Chapter 8]{Remy}
let 
$$\Delta_1:=\Delta^m(\{x,y\})\cup\Delta^u(\{x,y\})\cup -\Delta^u(\{x,y\})$$ and $$\Delta_2:=\Delta^m(\{x',y\})\cup\Delta^u(\{x',y\})\cup-\Delta^u(\{x',y\}).$$ Then $\Delta_1 \subseteq \Delta_2$ as $I \subseteq I'.$\\
It follows that $W_{x,y,\mA} \subseteq W_{x',y,\mA}.$ Since $H$ is contained in $\GL(W_{x,y,\mA})\leq \GL(W_{x',y,\mA})$, it follows that the Zariski closure of $H$, when computed in $\GL(W_{x,y,\mA})$ is the same as when computed in $\GL(W_{x',y,\mA})$, from which the claim follows.
\end{proof}

The independence from the fixed points $x,y$ allows us to associate a canonical subgroup to a bounded subgroup with trivial center. 
\begin{definition} Let $k$ be a field and let $G:=\GG_\DD(k)$ be a split Kac-Moody group. Let $U\leq G$ be a bounded subgroup with trivial center. Let $x\in \Delta_+,y \in \Delta_-$ be two points fixed by $U$. Then $\overline{\Ad U}:=U_{x,y}$ is called the {\bf Zariski closure of $U$.} \\
Let $\Addim U:=\dim \overline{\Ad U}$ denote the {\bf Ad-dimension} of $U$.
\end{definition}

We apply this to the study of the fixed point set.
\begin{lemma}\label{Addimfixpoint}
 Let $U,U'$ be two bounded subgroups with trivial center such that $U' \leq U.$ If $\Fix U \subsetneq \Fix U'$ then either there exists a finite index subgroup $U^*\leq U$ such that $\Fix U' \subseteq \Fix U^*$ or $\Addim U' < \Addim U.$
\end{lemma}
\begin{proof}
Suppose without loss of generality that $\Fix_U(\Delta_+) \subsetneq \Fix_{U'}(\Delta_+).$ Consider the CAT(0) realization of $\Delta_+$ and let 
 $x \in \Fix_{U'}(\Delta_+) \backslash \Fix_{U}(\Delta_+)$ and $y \in \Fix_U(\Delta_+).$ Consider a geodesic segment $p: [0,r] \to \Delta_+$ such that $p(0)=x, p(r)=y$. 
Let $s \in [0,r]$ be minimal such that $p(s) \in \Fix_{U}(\Delta_+)$. Let $z \in \Fix_U(\Delta_-)$ and let $P:=\Fix_G\{p(s), z\}.$ Then $U$ is contained in the bounded subgroup $P.$\\
Note that $s>0$ and that for some $\varepsilon>0$ the segment $[p(s-\varepsilon), p(s)]$ is contained in a residue $R'$, which is a proper residue of the spherical building associated to $\{p(s),z\}.$ This in
turn says that $\overline{\Ad U'}$ is contained in a proper parabolic $P'$ of $P$, while $\overline{\Ad U}$ is not. If the connected component of the identity of $\overline{\Ad U}$ is not contained in $P$, we must have $\Addim U>\Addim U'$ by \cite[1.8.2]{Springer} applied to  $\overline{\Ad U'}^0$ and $\overline{\Ad U}^0$.  \\
Otherwise let $U^*$ be the preimage of $\Ad U \cap \overline{\Ad U}^0$ in $U$. Then $U^*$ is of finite index in $U$ and $\overline{\Ad U^*}$ is connected. By the same reasoning, using the previous proposition, the fixed point set of the group $U^*$ must necessarily contain the fixed point set of $U'$.     
\end{proof}

\begin{remark}
Clearly, the fixed point set of a finite index subgroup $U'$ of $U$ can be much larger. For example, $\Sym(n)$ operates on $V:=k^n$ via permutation of the basis vectors and leaves the subspace generated by $v=(1, \ldots, 1)$ invariant. The induced action on $V/\langle v \rangle \cong k^{n-1}$ then is irreducible. It follows that $\Sym(n)$ is a subgroup of $\SL_{n-1}(k)$ which does not fix a proper residue of the spherical building associated to $\SL_{n-1}(k)$, while it virtually acts trivially.   
\end{remark}

\begin{theorem}\label{infinitealgebraicbounded} Let $L$ be an algebraic extension of $\QQ$. Let $k'$ be a field of characteristic 0 and $G:=\GG_\DD(k)$ a split Kac-Moody group. Then any homomorphism $\varphi\colon \SL_2(L) \to G(k)$ has bounded image. 
\end{theorem}

\begin{proof}
 If $L$ is a finite extension, the result follows from Corollary \ref{numberfieldbounded}. Otherwise let $(L_i)_{i \in \NN} \subseteq L$ be an ascending sequence of subfields such that $L=\bigcup L_i$. For each $i$, choose a set of valuations $S_i$ of $L_i$ such that $W_i:=\SL_2(L_i(S_i))$ is an irreducible lattice and such that $S_i \subseteq S_{i+1}.$ Let $U_i \leq W_i$ be a finite index subgroup with trivial center, as provided by Lemma \ref{Gamma0}.\\

For each $i$, $U_i(\leq \SL_2(L_i))$ is bounded by Corollary \ref{numberfieldbounded}. By Theorem \ref{MargulisThm}, the Zariski closure of $U_i$ is semisimple, hence by Lemma \ref{Hochschild} contained in a Levi factor. In particular, the Ad-dimension of $U_i$ is bounded above since there are only finitely many conjugacy classes of Levi factors in $G(K).$ \\
 Pick $i_0$ such that $\Addim U_{i_0}=\max \{\Addim U_i: i\in \NN\}$ and let $x$ be a fixed point of $U_{i_0}.$ 
Let $U_i^*$ be the finite index subgroup of $U_i$ provided by Proposition \ref{Addimfixpoint} which fixes $x.$ Then $U^*:=\langle U_i^*: i \in \NN \rangle$ fixes $x$ and for each $i$, the index of $U^* \cap U_i$ is finite. Let $\OO_i$ denote the ring of integers of $L_i$. Then $\SL_2(\OO_i)\leq W_i$ and hence $V_i:=U^* \cap \SL_2(\OO_i)$ has finite index in $\SL_2(\OO_i).$\\ 
We will show that $U^*$ and $\SL_2(\QQ)$ boundedly generate $\SL_2(L)$, which will imply the claim by Lemma \ref{boundedgenerationlemma}. Since every $g\in \SL_2(L)$ is contained in some $\SL_2(L_i)$, it suffices to prove that $\SL_2(\QQ)$ and $V_i$ boundedly generate $\SL_2(L_i)$ and uniformly so. This is the content of the following lemma. 
\end{proof}

Let $N_0$ denote the maximum number of elementary matrices needed to express a matrix $g\in G=\SL_2(K)$ as a product of elementary matrices. Note that each torus element $t$ is the product of at most 6 elementary matrices, e.g. $t=m(u)m(1)$. By the Bruhat decomposition, $G=TU_+ \cup U_+TsU_+$, so $N_0 \leq 11.$
\begin{lemma} Let $L$ be a number field and $\OO$ its ring of integers. Let $V$ be a subgroup of $\SL_2(\OO_i)$ of finite index. Then every element of $\SL_2(L)$ can be written as a product of at most $3N_0$ matrices from either $\SL_2(\QQ)$ or $V$. 
\end{lemma}

\begin{proof}
For $x$ in $L$ there is some $q \in \NN$ such that $qx \in \OO$. Since $\OO$ is a ring containing $\ZZ$, $q^2 x\in \OO$. Since $V$ has finite index in $\SL_2(\OO)$ there is some $a \in \NN$ such that 
$u_+(a^2q^2x) \in V.$
We may write $u_+(x)=\diag((aq)^{-1},aq)u_+(a^2q^2x)\diag(aq,(aq)^{-1})$, from which the claim follows.
\end{proof}

\bibliographystyle{plain}
\bibliography{References}
\end{document}